\documentclass[11pt, sec]{article}
\usepackage{graphicx}
\usepackage{amsmath}
\usepackage{amssymb}
\usepackage{amsfonts}
\usepackage{amscd}
\usepackage{latexsym}
\usepackage{graphics}

  \newtheorem{theorem}{Theorem}[section]
  \newtheorem{prop}[theorem]{Proposition}
  \newtheorem{cor}[theorem]{Corollary}
  \newtheorem{lemma}[theorem]{Lemma}
  \newtheorem{Remark}[theorem]{Remark}
  \newtheorem{Def}[theorem]{Definition}
  
  \newtheorem{thm}[theorem]{Theorem}

\def\deg{\operatorname{deg}}
\def\det{\operatorname{det}}

\def\Hom{\operatorname{Hom}}

\def\deg{\operatorname{deg}}

\newcommand{\lra}{\longrightarrow}

\newcommand{\mt}{\mapsto}
\newcommand{\h}{\simeq}

\newcommand{\mb}{\mathbb}

\newcommand{\f}{\frac}
\renewcommand{\l}{\left}
\renewcommand{\r}{\right}
\renewcommand{\lg}{\langle}
\newcommand{\rg}{\rangle}
\newcommand{\be}{\begin{equation}}
\newcommand{\ee}{\end{equation}}
\newcommand{\bce}{\begin{center}}
\newcommand{\ece}{\end{center}}
\newcommand{\bq}{\begin{eqnarray}}
\newcommand{\eq}{\end{eqnarray}}

\newcommand{\ra}{\rightarrow}

\renewcommand{\a}{\alpha}

\renewcommand{\rm}{\textrm}
\newcommand{\q}{\quad}
\newcommand{\qq}{\qquad}

\newcommand{\e}{\emph}
\newcommand{\la}{\lambda}

\newcommand{\s}{\sum}
\renewcommand{\i}{\infty}
\renewcommand{\c}{\cdot}
\newcommand{\cs}{\cdots}
\newcommand{\ba}{\begin{array}}
\newcommand{\ea}{\end{array}}

\newcommand{\om}{\omega}

\renewcommand{\e}{\emph}
\renewcommand{\s}{\sigma}

\renewcommand{\ss}{\subset}
\newcommand{\bpm}{\begin{pmatrix}}
\newcommand{\epm}{\end{pmatrix}}
\newcommand{\bbm}{\begin{bmatrix}}
\newcommand{\ebm}{\end{bmatrix}}
\newcommand{\ld}{\ldots}

\newcommand{\op}{\oplus}

\newcommand{\ot}{\otimes}
\newcommand{\hr}{\hookrightarrow}

\newcommand{\bop}{\bigoplus}

\newcommand{\wtl}{\widetilde}

\begin{document}

\title{\Large \bf $q$-Conjugacy Classes in Loop Groups     }
\author{Dongwen Liu}
\date{}

\maketitle

\section{\large Introduction}
Let $k((t))$ be the field of formal Laurent power series over some field $k$ of characteristic zero, $q\in k^\times $ which is not a root of unity, and $\s_q$ be the automorphism of $k((t))$ defined by
$\s_q: \a(t)\mt\a(qt).$ For an algebraic group $G$ over $k$ write
$G(k((t)))$ for the $k((t))$-points of $G$. Then $\s_q$ also acts on $G(k((t)))$ by change of variable, namely $\s_q(g(t))= g(qt).$ One can define a $q$-conjugate action of $G(k((t)))$ on itself by
\[
\rm{Ad}_qh:\q g\longmapsto h g\s_q (h)^{-1}.
\]

The main purpose of this paper is to classify the $q$-conjugacy classes in $G(k((t)))$ for $G=GL_n,$ $SL_n,$ $Sp_{2n}$ or $O_n.$ Our main tool is the so-called $D_q$-modules, and the classification is based on the principle that there is a one-to-one correspondence between the set of isomorphism classes of $n$-dimensional $D_q$-modules and the set of $q$-conjugacy classes in $GL_n(k((t)))$. Similarly if we equip $D_q$-modules with symplectic or symmetric forms, then we are also able to classify $q$-conjugacy classes in $Sp_{2n}(k((t)))$ and $O_n(k((t)))$.

The motivation to classify the $q$-conjugacy classes in loop groups is that we can interpret them as ordinary conjugacy classes in a larger group. As explained in \cite{BG}, consider the semidirect product $G(k((t)))\rtimes k^\times $ where $q\in k^\times$ acts on $G(k((t)))$ as $\s_q.$ Note that $x, y\in G(k((t)))$ are $q$-conjugate if and only if the elements $x\s_q$ and $y\s_q$ of $G(k((t)))\rtimes \lg \s_q\rg$ are conjugate under $G(k((t)))$. A detailed study of conjugacy classes in $G(k((t)))\rtimes k^\times $, as well as its central extensions, ought to be quite useful for generalizing the theory of automorphic forms to infinite dimensional groups.

Another important feature of $q$-conjugacy classes is their geometric interpretation. More precisely, there exist very close connections between $q$-conjugacy classes in $G(\mb{C}((t)))$ and holomorphic principle $G$-bundles on Tate's elliptic curve $E_q=\mb{G}_m/ q^\mb{Z}$ over $\mb{C}$, where $q\in \mb{C}^\times $ and $|q|<1.$ $G$-bundles over elliptic curves have been well understood, see \cite{A}, \cite{FM}. We shall summarize some results along this line in the last section.

This paper is organized as follows. In section 2 we first review the results on $D_q$-modules in \cite{vdP-R} which essentially classify the $q$-conjugacy classes in $GL_n(\mb{C}((t))).$ We shall generalize the results to the case that the ground field is not algebraically closed, using some Galois arguments. We also discuss the classification in the central extension of $GL_n(k((t))$, and introduce a new category of $D_q$-modules. Section 3 deals with symplectic and orthogonal groups and gives the classification theorem. In the last section we discuss principal bundles over elliptic curves and we shall focus on regular semi-stable case.
\\\\
{\bf Acknowledgements.} The author thanks his advisor, Professor Y.-C. Zhu for suggesting this question, and for many helpful discussions and  encouragement during this work.

\section{\large $D_q$-modules and $q$-conjugacy classes in $GL_n$}

\subsection{Loop groups without central extensions}

First assume that $k$ is algebraically closed of characteristic zero and $q\in k^\times$ is not a root of unity. Let us recall the definition and classification of
 $D_q$-modules or the so-called difference modules, following \cite{BG}, \cite{vdP-R}, \cite{Sau1}, \cite{Sau2}. Equivalently this classifies the
 $q$-conjugacy classes in $GL_n(k((t)))$ in this case. Later we will consider the non-algebraically closed case.

\begin{Def}
A $D_q$-module over $k((t))$ is a pair $(V, \Phi)$, where $V$ is a finite dimensional vector space over $k((t))$ and $\Phi$ is a $\s_q$-semilinear bijection $V\ra V$, i.e.
\[
\Phi(\a v)=\s_q (\a)\Phi(v),\q \forall\a\in k((t)), v\in V,
\]
where $\s_q\in \mathrm{Aut}_kk((t))$ is given by $t\mt q t.$
\end{Def}

$D_q$-modules over $k((t))$ form a tensor category ${\cal D}_q(k,t)$, where morphisms between two objects $(V,\Phi)$ and $(V',\Phi')$ are all the $k((t))$-linear maps $\rho: V\ra V'$ such that $\Phi'\circ\rho=\rho\circ\Phi.$

\begin{prop}\label{1}
There is a bijection between the set of $q$-conjugacy classes in $GL_n(k((t)))$ and the set of isomorphism classes of $n$-dimensional $D_q$-modules over $k((t))$.
\end{prop}
\e{Proof.} Let $V$ be an $n$-dimensional vector space over $k$ and $V_{k((t))}=V\ot_k k((t)).$ For any element $g\in GL(V_{k((t))})$ we get a $D_q$-module $(V_{k((t))},\Phi_g),$ where $\Phi_g=g\circ (\rm{id}_V\ot \s_q).$ If $g'$ is $q$-conjugate to $g$, i.e. $g'=h g \s_q (h)^{-1}$ for some $h\in G(k((t)))$, then the two $D_q$-modules we get are isomorphic. Conversely, given a $D_q$-module $(W,\Phi)$, the matrix of $\Phi$ under a basis of $W$ gives an element $g$ in $GL_n(k((t)))$. Another choice of the basis gives a $q$-conjugate of $g.$ Then we get a $q$-conjugacy class in $GL_n(k((t)))$. Two $D_q$-modules give the same class if and only if they are isomorphic. It is easy to check that the two constructions are inverse to each other and give the required bijection.
\hfill$\Box$

\

The following cyclic vector lemma is well-known, which essentially appears in Birkhoff's work \cite{B}. \cite{DV} includes an elegant proof. See also
 \cite{K}.

\begin{lemma}
Let $(V,\Phi)$ be an $n$-dimensional $D_q$-module over $k((t))$, then $V$ contains a
cyclic vector $v$, i.e. $\{v,\Phi v,\ld,\Phi^{n-1}v\}$ is a
$k((t))$-basis of $V$.
\end{lemma}

Let $k((t))_q[T]$ be the non-commutative polynomial ring in one indeterminate $T$ over $k((t))$, whose underlying abelian group is the same as $k((t))[T]$ but the multiplication rule is given by
$T \a =\s_q(\a) T$ for all $\a\in k((t)).$ The lemma implies that $V$ is isomorphic to $k((t))_q[T]/(f(T))$ for some $f(T)=T^n+a_1T^{n-1}+\cs+a_n$ with
 $a_i\in k((t))$ and $a_n\neq 0$, and $\Phi$ is given by left multiplication by $T.$

\

Before we proceed we would like to digress a bit and review some other familiar examples.

(i) If $V$ is a finite dimensional vector spaces over some field $k$ and $\Phi: V\ra V$ is in $\mathrm{End}_k(V)$, then $V$ under $\Phi$ decomposes into a direct sum of submodules, each of which is isomorphic to $k[T]/(f^e)$ for some irreducible polynomial $f(T)\in k[T]$ and some integer $e\geq 1,$ and $\Phi$ acts by multiplication by $T$. If we assume that $\Phi\in \mathrm{Aut}_k(V),$ then $f(T)\neq T$ in each summand of the decomposition. The theory of Jordan canonical forms classifies ordinary conjugacy classes in $GL(V).$

(ii) Let $k$ be a finite field, $\s\in\rm{Gal}(\bar{k}/k)$ be the Frobenius element. $g$ and $g'$ in $GL_n(\bar{k})$ are said to be $\s$-conjugate if $g'=h g \s(h)^{-1}$ for some $h\in GL_n(\bar{k}).$ Then Lang's lemma says that there is only one $\s$-conjugacy class in $GL_n(\bar{k}).$

(iii) Let $k$ be a perfect field of characteristic $p>0,$ $W=W(k)$ the ring of Witt vectors over $k$, $\s: W\ra W$ the Frobenius endomorphism of $W$, $K$ the fraction field of $W.$ An $F$-isocrystal on $k$ is a finite dimensional vector space $V$ over $K$, together with a $\s$-semilinear bijection $\Phi: V\ra V.$ Let $K_\s[T]$ be the non-commutative polynomial ring in one indeterminate $T$ over $K$, with multiplication rule $T \a=\s(\a) T$ for all $\a\in K.$ Following the proof of Proposition \ref{1}, we can construct a bijection from the set of $\s$-conjugacy classes in $GL_n(K)$ to the set of isomorphism classes of $n$-dimensional $F$-isocrystals.

Dieudonn\'{e} \cite{D} and Manin \cite{M} proved that if $k$ is algebraically closed, then the category of $F$-isocrystals is semi-simple and its simple objects are $V_{r,s}=K_\s[T]/(T^s-p^r),$ where $r,s\in\mb{Z},$ $s\geq 1,$ $(r,s)=1,$ and $\Phi$ is given by left multiplication by $T.$

The rational number $\a=r/s$ is called the \e{slope} of $V_{r,s}.$ Any $F$-isocrystal $V$ splits into a direct sum
$V=\bop\limits_{\a\in\mb{Q}}V_\a$, where $V_\a$ is a direct sum of simple $F$-isocrystals of slope $\a.$ The dimension of $V_\a$ over $K$ is called the \e{multiplicity} of $\a.$ If $\a_1<\a_2<\cs<\a_k$ are slopes of $V$ and $\la_i$ is the multiplicity of the slope $\a_i$, then the \e{Newton polygon} of $V$ is the line joining the points $P_0=(0,0),$ $P_i=(\sum^i_{j=1}\la_i,\sum^i_{j=1}\la_j\a_j),$ $i=1,\ld, k.$ It is known that if $k=\mb{F}_{p^m}$ then the Newton polygon of $V$ is the ordinary Newton polygon of the polynomial $\det(t\c\rm{Id}-\Phi^m).$

\

In the same spirit as the last example, the notion of slopes and Newton polygons were also defined for $D_q$-modules and we shall describe them below. The comparisons between $F$-isocrystals and $D_q$-modules are quite instrumental.

Let $\nu: k((t))\ra\mb{Z}\cup\{\i\}$ be the valuation with respect to $t$. For any $D_q$-module $V$, by the cyclic vector lemma we can assume that it is of the form $k((t))_q[T]/(f(T))$, where $f(T)=T^n+a_1T^{n-1}+\cs+a_n$ with $a_n\neq0$. Consider the convex hull of $\bigcup^n_{i=0}\{(i,\nu(a_i)+x)|x\leq 0\}.$ The finite  part of the boundary of this hull is called the \e{Newton polygon} of $V.$ $V$ is called \e{pure} if its Newton polygon has only one slope. Any $D_q$-module splits uniquely into direct sum of pure modules. In particular any irreducible $D_q$-module is pure. See \cite{Sau1}, \cite{Sau2} for details.

From a geometric point of view, one can identify the notion of pure with \e{semi-stable}, and construct the Newton polygon as follows. For $V\in {\cal D}_q(k,t)$ pick up an element $g$ in the corresponding conjugacy class in $GL_n(k((t)))$ where $n=\rm{rank }V.$ Define the degree of $V$ to be $\deg V:=\nu(\det g).$ Clearly the definition does not depend on the choice of $g.$
  For $V\neq 0$ introduce
 the rational number
\[
\mu(V)=\deg V/ \rm{dim}V.
\]
$V$ is called semi-stable (resp. stable) if for any nonzero proper $D_q$-submodule $V'$ we have $\mu(V')\leq (\rm{ resp.} <)\mu(V)$. Consider the
 convex hull of $\{(\rm{dim} V', \deg V'+x)| x\leq 0, V'\subseteq V\rm{ a }D_q\rm{-submodule}\}.$ Then the Newton polygon is the finite part of the boundary of
 the convex hull.

This gives a canonical filtration, which is analogous to the well-known Harder-Narasimhan filtration \cite{HN}  \cite{Sau2}
\[
0=V_0\ss V_1\ss \cs \ss V_l=V
\]
such that

(i) $V_i/V_{i-1}$ is semi-stable for $i=1,\ld, l$,

(ii) $\mu(V_i/V_{i-1})>\mu(V_{i+1}/V_i)$ for $i=1,\ld, l-1.$
\\
The graded module $gr V$ associated to $V$ is $\bigoplus^l_{i=1}V_i/V_{i-1}.$ It is known that $V\h gr V$. The notion of stability coincides with that of vector bundles on elliptic curves, which will be discussed in the last section.

\

Let us describe the classification of $D_q$-modules \cite{vdP-R}.
Fix a compatible system of roots of $q$, i.e. $\{q^{1/n}|n\in
\mb{Z}\}$ such that $(q^{1/n})^m=q^{m/n}$ whenever $m|n.$ Then we
can extend the action of $\s_q$ to $k((t^{1/n}))$ by $\s_q
t^{1/n}=q^{1/n}t^{1/n}.$ Then scalar extension gives a functor
\[
ext_n:\q {\cal D}_q(k,t)\longrightarrow{\cal D}_{q^{1/n}}(k,t^{1/n}).
\]
In fact we may extend the action of $\s_q$ to $\overline{k((t))}=\bigcup_{n\geq 1}k((t^{1/n})$ and obtain an obvious functor
\[
ext_\i:\q {\cal D}_q(k,t)\longrightarrow\lim_n{\cal D}_{q^{1/n}}(k,t^{1/n}).
\]
Conversely restriction of scalar from $k((t^{1/n}))$ to $k((t))$ also gives a functor
\[
res_n:\q {\cal D}_{q^{1/n}}(k,t^{1/n})\longrightarrow{\cal D}_q(k,t).
\]
The building blocks of $D_q$-modules are the following objects.

Consider $r,s\in\mb{Z}$ with $s\geq 1$ and $(r,s)=1$. For any
$\bar{\la}\in k^\times/q^\mb{Z}$ pick up a representative $\la\in
k^\times.$ Let $E(\la^{1/s} t^{r/s})=k((t^{1/s}))_q[T]/(T-\la^{1/s}t^{r/s})\in{\cal
D}_{q^{1/s}}(k,t^{1/s})$, and let $W_{r,s}(\la)=res_s(E(\la^{1/s} t^{r/s})).$ Equivalently, $W_{r,s}(\la)$ is isomorphic to $k((t))_q[T]/(T^s-q^{r(s-1)/2}\la t^r).$

For convenience we denote  by $V_{r,s}(\la)\in{\cal D}_q(k,t)$ the $D_q$-module $k((t))_q[T]/(T^s-\la t^r).$
Let $J_d$ be the $D_q$-module $k((t))_q[T]/((T-1)^d),$ which is the analogue of unipotent Jordan blocks.

\begin{thm}\label{PR} $\mathrm{(M.~van~der~Put,~M.~Reversat~[12])}$

$\mathrm{(i)}$ Up to isomorphism $V_{r,s}(\la)$ only depends on $r, s$ and $\bar{\la}$, and is irreducible with slope $r/s.$ This gives the isomorphism
classes of irreducible $D_q$-modules over $k((t)).$

$\mathrm{(ii)}$ Any indecomposable $D_q$-module over $k((t))$ is isomorphic to some $J_d\otimes V_{r,s}(\la)$ with
$(d,r,s,\bar{\la})$ uniquely determined.
\end{thm}

It is clear from the definition and Theorem \ref{PR} that stable $D_q$-modules are irreducible ones, and semi-stable $D_q$-modules with $\mu=r/s$
 are of the form $\bigoplus^l_{i=1}J_{d_i}\otimes V_{r,s}(\la_i).$ Notice that the category ${\cal D}_q(k,t)$ is not semi-simple, different from the category of $F$-isocrystals.

\begin{cor}\label{cor}
$q$-conjugacy classes in $GL_n(k((t)))$ are parameterized by isomorphism classes of $n$-dimensional $D_q$-modules over $k((t))$ of the form $\bigoplus^l_{i=1}J_{d_i}\otimes V_{r_i,s_i}(\la_i).$
\end{cor}

\begin{prop}
The natural map
\[
q\rm{-conjugacy classes in }SL_n(k((t)))\lra q\rm{-conjugacy classes in }GL_n(k((t)))
\]
is injective.
\end{prop}
\e{Proof.} Suppose $a,b\in SL_n(k((t)))$, $c\in GL_n(k((t)))$ such that $b= c\c a\c \s_q c^{-1}.$ Taking determinants, we see that $\det c$ is fixed by
 $\s_q$, i.e. $\det c\in k.$ Set $c'=(\det c)^{-1/n}c\in SL(n,k((t)))$, then $b= c'\c a\c \s_q c'^{-1}$.
\hfill $\Box$

\begin{cor}
$q$-conjugacy classes in $SL_n(k((t)))$ are parameterized by classes of $D_q$-modules of rank $n$ of the form $\bigoplus^l_{i=1}J_{d_i}\otimes V_{r_i,s_i}(\la_i)$ such that $\sum_i d_i r_i=0$ and $\prod_i \la_i^{d_i}\in q^\mb{Z}.$
\end{cor}
\e{Proof.} $g\in GL_n(k((t)))$ can be conjugated into $SL_n(k((t)))$ if and only if $\det g$ can be written in the form  $q^ru,$ where $r\in\mb{Z},$ $u=1+a_1 t+a_2t^2+\cs, a_i\in k.$ It suffices to consider $GL_1$. Indeed, there exists $v\in k[[t]]^\times$ such that $u=v(\s_q v)^{-1}.$ Then $q^r u= (t^{-r} v)\s_q(t^{-r}v)^{-1}.$
\hfill $\Box$

\

In the rest of this subsection, we shall only assume that $k$ is of characteristic zero, which is not necessarily algebraically closed, and $q\in k^\times$ again is not a root of unity. Note that the cyclic lemma holds for
any $k$ which is infinite.

We shall use some simple Galois descent argument. Fix an algebraic closure $\bar{k}$ of $k$. Note that $\bar{k}((t))$ is not algebraic over $k((t))$, although ${\cal D}_q(\bar{k},t)$ has been
classified. Let us introduce $\bar{k}((t))_0:=k((t))\ot_k \bar{k}$, which is algebraic over $k((t))$ and $\rm{Gal}(\bar{k}((t))_0/k((t)))\h \rm{Gal}(\bar{k}/k).$

\begin{lemma}
The $D_q$-modules over $\bar{k}((t))_0$ have the same classification as Theorem \ref{PR}.
\end{lemma}
\e{Proof.} The argument in \cite{vdP-R} still works, once we prove that any one-dimensional $D_q$-module over
$\bar{k}((t))_0$ is isomorphic to $W_1(r,\la)$ for some $r\in\mb{Z}$ and $\la\in k^\times.$ This amounts to the fact that the multiplicative group
$U:=\{1+a_1t+a_2t^2+\cs|a_i\in \bar{k}\}\cap \bar{k}((t))_0$ enjoys the property that any $u\in U$ can be written as $v^{-1}\s_q v$ for some $v\in U.$
\hfill $\Box$

\

Let ${\cal D}_q(\bar{k},t)_0$ be the category of $D_q$-modules over
$\bar{k}((t))_0$ and $ext_{\bar{k}}$ be the functor of scalar
extension:
\[
ext_{\bar{k}}:\q {\cal D}_q(k,t) \stackrel{\ot_k \bar{k}}{\lra}
{\cal D}_q(\bar{k},t)_0.
\] Again we denote by $V_{r,s}(\la)\in{\cal D}_q(\bar{k},t)_0$ the $D_q$-module
$\bar{k}((t))_{0,q}[T]/(T^s-\la t^r)$, $\la\in\bar{k}^\times,$ where $\bar{k}((t))_{0,q}[T]=k((t))_q[T]\ot_k \bar{k}$,
$r,s\in\mb{Z},$ $(r,s)=1,$ $s\geq 1$. Let $S$ be the set of monic irreducible polynomials over $k$ with non-zero constant terms. Then we have the bijection \[S\stackrel{\sim}{\lra} \bar{k}^\times/\rm{Gal}(\bar{k}/k),\]
where $\bar{k}^\times/\rm{Gal}(\bar{k}/k)$ is the set of Galois orbits in $\bar{k}^\times$ under the action of $\rm{Gal}(\bar{k}/k).$
Indeed associated to $f\in S$ there is a Galois orbit ${\cal O}_f$ in $\bar{k}^\times$, which is the set of roots of $f$ in $\bar{k}^\times.$
There is an action of $\s_q$ on $S$, defined by $\s_q f(x)= q^{-\deg f}f(qx).$ Corresponding to the bijection above, $\s_q$ acts on $\bar{k}^\times/\rm{Gal}(\bar{k}/k)$ by $\s_q{\cal O}=q^{-1}{\cal O}.$ Modulo the action of $\s_q$ we get
\[
S/\lg \s_q\rg \stackrel{\sim}{\lra}
(\bar{k}^\times/\rm{Gal}(\bar{k}/k))/\lg \s_q\rg.
\]
For $\bar{f}\in S/\lg \s_q\rg$ let $f(x)\in S$ be a
representative. For $r,s\in\mb{Z},$ $(r,s)=1,$ $s\geq 1$ let $V_{r,s}(f)\in{\cal
D}_q(k,t)$ be the $D_q$-module $k((t))_q[T]/(f(t^{-r}T^s)).$ In particular, we have $V_{r,s}(\la)=V_{r,s}(x-\la)$ for $\la\in k.$

More explicitly, if we write $f(x)=x^m-a_{m-1}x^{m-1}-\cs-a_0$, then
\[
f(t^{-r}T^s)=q^{-\f{1}{2}srm(m-1)}t^{-rm}T^{sm}-\sum^{m-1}_{j=0}q^{-\f{1}{2}srj(j-1)}a_j t^{-rj}T^{sj}.
\]
$V_{r,s}(f)$ can be identified with the $D_q$-module
with basis $\{\Phi^{i}v_j|i=0,\ld,s-1,j=0,\ld,m-1\}$ such that
\[
\Phi^s v_j=t^r v_{j+1}, \q j=0,\ld, m-2,\q \Phi^s v_{m-1}=t^r(a_0
v_0+\cs+a_{m-1}v_{m-1}).
\]
The equivalence of the two presentations is given by $v_0\mt T.$ Note that $v_0$ is a cyclic vector.

\

Our classification is the following

\begin{thm}\label{des}

$\mathrm{(i)}$ Up to isomorphism $V_{r,s}(f)$ only depends on $r, s$ and
$\bar{f}\in S/\lg\s_q\rg,$ and is irreducible. This gives the isomorphism
classes of irreducible $D_q$-modules over $k((t)).$

$\mathrm{(ii)}$ $ext_{\bar{k}}(V_{r,s},f)\h \bigoplus_{\la\in{\cal O}_f}
V_{r,s}(\la),$ where ${\cal O}_f$ is the set of roots of $f.$

$\mathrm{(iii)}$ Any indecomposable $D_q$-modules over
$k((t))$ is isomorphic to some $J_d\ot V_{r,s}(f)$ with $(d, r, s, \bar{f})$
uniquely determined.
\end{thm}
\e{Proof.} (ii) is easy. Use above description of $V_{r,s}(f)$, in the $\bar{k}$-span of $\{v_0,\ld,
v_{m-1}\}$ one can find another basis $\{w_1,\ld w_m\}$ such that
$\Phi^n w_i=\la_it^r w_i$ and ${\cal O}_f=\{\la_1,\ld,\la_m\}.$
The proof of (i) consists of two steps.

(a) $V_{r,s}(f)$ is irreducible: Let $U$ be a nonzero $D_q$-submodule of $V_{r,s}(f).$ With notation as above, $\rm{Gal}(\bar{k}/k)$ acts on the $\bar{k}$-span of $\{v_0,\ld,v_{m-1}\}$, and we may assume that $\rm{Gal}(\bar{k}/k)$ permutes $\{w_1,\ld,w_m\}$ transitively. This implies that $\rm{Gal}(\bar{k}/k)$ permutes the summands in (ii) transitively. Consider the composition
\[
\varphi_i:\q ext_{\bar{k}}(U)\hookrightarrow ext_{\bar{k}}(V_{r,s}(f))\stackrel{\pi_i}{\lra} V_{r,s}(\la_i).
\]
One of $\varphi_i$'s hence all of them are surjective. It follows easily that $ext_{\bar{k}}(U)= ext_{\bar{k}}(V_{r,s}(f))$ by considering the decomposition of
$ext_{\bar{k}}(U).$

(b) Any irreducible object $V\in {\cal D}_q(k,t)$ is of this form, and $V_{r,s}(f)\h V_{r,s}(g)$ if $\bar{f}=\bar{g}\in S/\lg \s_q\rg$: Take an irreducible submodule of $ext_{\bar{k}}(V)\in {\cal D}_q(\bar{k},t)_0,$ which is of the form $V_{r,s}(\la)$ for some $\la\in\bar{k}^\times.$ Let $f\in S$ be the minimal polynomial of $\la$ over $k$, and let $K$ be the Galois extension of $k$ generated by all roots of $f.$ Let $v$ be the cyclic vector of $V_{r,s}(\la)$ so that $\Phi^sv=\la t^rv.$ Then $v$ is defined over a finite extension $L$ of $K$, i.e. $v\in V\ot_k L.$ We may assume $L/K$ is Galois and then there exist some $\a \in L$ such that $v':=\sum\limits_{\s\in \rm{Gal}(L/K)}\s (\a v)\neq 0.$ Let us still write $v$ for $v'$, then $v\in V\ot_k K$ and $\Phi^s v=\la t^r v.$

Utilize the same trick again we may assume that $\sum\limits_{\s\in \rm{Gal}(K/k)} \s v\neq 0.$ Define
\[
u_j=\sum_{\s\in \rm{Gal}(K/k)}\s(\la^j v),\q j=0,\ld, m-1.
\]
Then $u_j\in V$ and satisfy the relations
\[
\Phi^s u_j=t^r u_{j+1}, \q j=0,\ld, m-2,\q \Phi^s u_{m-1}=t^r(a_0
u_{0}+\cs+a_{m-1}u_{m-1}).
\]
This implies that $\Phi^iu_j$'s span $V$ since $V$ is irreducible, and there is a nonzero morphism from $V_{r,s}(f)$ to $V$ which maps $v_j$ to $u_{j}.$ Therefore (ii) implies $V\h V_{r,s}(f).$ If $g\in S$ and $\bar{f}=\bar{g}\in S/\lg\s_q\rg$, then $\la\equiv \mu \mod q^\mb{Z}$ for some $\mu\in{\cal O}_g.$ Since $V_{r,s}(\la)\h V_{r,s}(\mu)$, above proof also implies that $V\h V_{r,s}(g).$ Then $V_{r,s}(f)\h V\h V_{r,s}(g).$

We only sketch the proof of (iii). Assume that $V\in{\cal D}_q(k,t)$ is indecomposable, and that $ext_{\bar{k}}(V)
\h \sum_i J_{d_i}\ot V_{r_i,s_i}(\la_i).$ The action of $\rm{Gal}(\bar{k}/k)$ permutes the quadruples $\{(d_i, r_i, s_i, \la_i)\}$, and indecomposability of $V$ implies that the action has a single orbit. Therefore we can assume that
$ext_{\bar{k}}(V)\h \sum_{\la\in{\cal O}_f}J_d\ot V_{r,s}(\la)$ for some $f\in S.$ Then use some simple calculations one can prove that $V\h J_d\ot V_{r,s}(f).$
\hfill $\Box$

\

Let us end this subsection with some remarks and computations on Galois cohomology. The following result follows from a standard argument \cite{Se}.

\begin{prop}\label{se}
Given $V\in {\cal D}_q (k,t)$, there is a bijective map from the isomorphism classes of $V'\in {\cal D}_q(k,t)$ such that $ext_{\bar{k}}(V')\h ext_{\bar{k}}(V)$ to $H^1 (\mathrm{Gal}(\bar{k}/k), \mathrm{Aut}( ext_{\bar{k}}(V))).$
\end{prop}

Recall that $\s\in \rm{Gal}(\bar{k}/k)$ acts on $\mathrm{Aut}(ext_{\bar{k}}(V)))$ by $\eta\mt \s(\eta)=\s\circ \eta\circ \s^{-1}.$ Theorem \ref{des} implies that  $H^1 (\rm{Gal}(\bar{k}/k), \mathrm{Aut}( ext_{\bar{k}}(V)))$ is trivial for any irreducible $V\in {\cal D}_q(k,t).$ It would be interesting to investigate the cocycles and coboundaries explicitly.
\begin{lemma}\label{aut}
$\mathrm{End}(V_{r,s}(\la))=\bar{k}$ for any $\la\in\bar{k}^\times.$
\end{lemma}
\e{Proof.} It is equivalent to show $\mathrm{Aut}(V_{r,s}(\la))=\bar{k}^\times.$ A direct proof is to compute the twisted centralizer in $GL_n(\bar{k}((t)))$ of the matrix representing $V_{r,s}(\la).$ Here we use another approach which reveals more results. Let
\[
U=\{u\in V_{r,s}(\la)| \Phi^s u=\la t^r u\}.
\]
Clearly $U$ contains a cyclic vector $v$, hence $\rm{End}(V_{r,s}(\la))\hr \rm{End}_{\bar{k}}(U)$ by restriction. We shall show that $U$ is one-dimensional over $\bar{k}$. In fact assume that
\[
u=a_0 v+a_1 \Phi v+\cs +a_{s-1}\Phi^{s-1}v\in U,\q a_i\in \bar{k}((t)).
\]
Then apply $\Phi^s$ to both sides we get
\[
\s_q^s a_i=q^{-ri}a_i,\q i=0,\ld,s-1.
\]
Since $(r,s)=1$ it follows that $a_i=0$ for $i>0$ and $a_0\in\bar{k}.$ Hence $U=\bar{k}v$.
\hfill $\Box$

\

For simplicity assume that $k(\la)$ is a Galois extension of $k$, and let $G=\rm{Gal}(k(\la)/k)$, then
\[
ext_{\bar{k}}(V)\h \bigoplus_{\s\in G}V_{r,s}(\s\la).
\]
The lemma implies that
\[
\rm{Aut}(ext_{\bar{k}}(V))=\bigoplus_{\s\in G}\bar{k}^\times.
\]
The action of $\rm{Gal}(\bar{k}/k)$ factors through $G$, and is given by
\[
\s:\q (f_\tau)_\tau\mt (\s f_{\s^{-1}\tau})_\tau,\q (f_\tau)_\tau\in \bigoplus_{\tau\in G}\bar{k}^\times.
\]
To compute cocycles and coboundaries, for a map $f: G\ra \rm{Aut}(ext_{\bar{k}}(V))$ let us write
$f(\s)=(f^\s_\tau)_\tau.$ Then $f$ is a coboundary if there exists some $(g_\tau)_\tau$ such that
\[
f_\tau^\s=\s g_{\s^{-1}\tau}\c g_\tau^{-1}.
\]
Some simple calculation shows that if $f$ is a cocycle then
\[
f^\s_\tau=\tau\l(f^{\tau^{-1}\s}_1(f_1^{\tau^{-1}})^{-1}\r).
\]
We can choose $g_\tau=\tau f^{\tau^{-1}}_1$. Therefore every cocycle is a coboundary, and we verified the triviality of
$H^1 (\rm{Gal}(\bar{k}/k), \mathrm{Aut}( ext_{\bar{k}}(V)))$ in this case.

\subsection{Central extensions of loop groups}

It is known that in general a Chevalley group $G$ over $k((t))$ has nontrivial central extensions
\[
1\lra k^\times \lra \wtl{G}(k((t))) \stackrel{\pi}{\lra} G(k((t)))\lra 1,
\]
e.g. via the tame symbols. We shall not discuss the general case but use the following construction for
$G=GL_n.$

A \e{lattice} $L$ of a $n$-dimensional $k((t))$-vector space $V$ is a free $k[[t]]$-submodule of rank $n$. In other words $L$ is a $k[[t]]$-span of a basis of $V$. Any two lattices $L_1$, $L_2$ in $V$ are \e{commensurable}, which means the quotients $L_1/(L_1\cap L_2)$ and $L_2/(L_1\cap L_2)$ are finite-dimensional over $k$. For example any lattice in $k((t))^n$ is commensurable with $k[[t]]^n.$

Let $L_0$ be the lattice $k[[t]]^n.$ Since $gL_0/(L_0\cap gL_0)$ is finite-dimensional over $k$, we can define the top wedge power $\wedge^{top} (gL_0/L_0\cap gL_0)$, which is a one-dimensional vector space over $k.$ Let $\det(L_0,gL_0)$ be the tensor product
\[
\wedge^{top} (gL_0/L_0\cap gL_0)\ot_k \wedge^{top} (L_0/L_0\cap gL_0)^{-1},
\]
where $(~)^{-1}$ denotes the dual vector space. And let $\det(L_0,gL_0)^\times$ be the set of nonzero vectors in $\det(L_0,gL_0)$, which form a $k^\times$-torsor.
Now define the group
\[
\wtl{GL}_n(k((t)))=\{(g, \om_g): g\in GL_n(k((t))), \om_g\in \det(L_0, gL_0)^\times\}.
\]
The multiplication in the group is given by
\[
(g,\om_g)(h,\om_h)=(gh, \om_g\wedge g\om_h),
\]
where $g\om_h$ is image of $\om_h$ under the natural map $\det(L_0, hL_0)\stackrel{g}{\ra} \det (gL_0, ghL_0)$, and
$\om_g\wedge g\om_h$ is defined by the isomorphism $\det(L_0,gL_0)\wedge \det(gL_0, ghL_0)\ra \det (L_0, ghL_0).$

$\s_q$ preserves $L_0$, hence induces the maps
  \[
  gL_0/L_0\cap gL_0\lra (\s_qg)L_0/L_0\cap (\s_qg)L_0,\q L_0/L_0\cap gL_0\lra L_0/L_0\cap (\s_qg)L_0.
  \]
  Therefore $\s_q$ induces the map
  $\det(L_0,gL_0)\ra \det(L_0,(\s_q g)L_0)$,
 hence acts on the group $\wtl{GL}_n(k((t))).$ To describe the $q$-conjugacy classes in $\wtl{GL}_n(k((t)))$, we introduce the following category $\wtl{{\cal D}}_q(k,t)$.

An object of $\wtl{{\cal D}}_q(k,t)$ is a quadruple $(V,\Phi, L, \om)$ where $(V,\Phi)$ form a $D_q$-module, $L$ is a lattice of $V$ and $\om\in \det (L,\Phi L)^\times.$ A morphism between two objects $(V,\Phi, L, \om)$ and $(V',\Phi', L', \om')$ is given by a pair $(\rho, u)$ where $\rho: (V,\Phi)\ra (V,\Phi')$ is a morphism of $D_q$-modules and $u\in
\det (L',\rho L)$ such that $u\wedge \rho\om\wedge \Phi' u^{-1}=\om'.$

Let us define the composition of morphisms. If $(\rho',u'): (V',\Phi', L', \om')\ra (V'',\Phi'', L'', \om'')$ is another morphism in $\wtl{{\cal D}}_q(k,t)$, then the composition $(\rho'',u'')=(\rho',u')\circ(\rho,u)$ is given by the pair $(\rho'\circ\rho, u'\wedge \rho'u).$ It is easy to check that $u''\wedge \rho''\om\wedge \Phi'' u''^{-1}=\om''$.

\begin{prop}\label{ext}
There is a bijection between the set of $q$-conjugacy classes in $\wtl{GL}_n(k((t)))$ and the set of isomorphism classes of $n$-dimensional $D_q$-modules in $\wtl{{\cal D}}_q(k,t).$
\end{prop}
\e{Proof.} Let $V$ be an $n$-dimensional vector space over $k,$ and $V_{k((t))}=V\ot_k k((t)).$ Then $\wtl{GL}_n(k((t)))$ can be realized as the group $\wtl{GL}(V_{k((t))})$, which consists of all the pairs $(g,\om_g)$ with $g\in GL(V_{k((t))})$ and $\om_g\in \det(V_{k[[t]]},gV_{k[[t]]})^\times$, where $V_{k[[t]]}:=V\ot_k k[[t]].$ Corresponding to $(g,\om_g)\in\wtl{GL}(V_{k((t))})$ we get a $D_q$-module $(V_{k((t))}, \Phi_g, V_{k[[t]]},\om_g)$ where $V_{k((t))}$ and $\Phi_g$ are defined as in the proof of Proposition \ref{1}. Notice that $gV_{k[[t]]}=\Phi_g V_{k[[t]]}$ since $\rm{id}_V\ot \s_q$ preserves $V_{k[[t]]}.$ It is easy to check that $D_q$-modules we get from $q$-conjugates of $(g,\om_g)$ are isomorphic.

Conversely, given an object $(W,\Phi,L,\om)$ in $\wtl{{\cal D}}_q(k,t),$ choose a basis of $L$ and use them to identify $L$ with $L_0=k[[t]]^n.$ Let $g$ be the matrix of $\Phi$ under this basis, and define $\om_g$ to be the image of $\om$ under the isomorphism $\det(L,\Phi L)\h \det(L_0,gL_0).$ Then different choices of the basis of $L$ give $q$-conjugates of $(g,\om_g)$ in $\wtl{GL}_n(k((t))).$ Two objects in $\wtl{{\cal D}}_q(k,t)$ gives the same $q$-conjugacy class if and only if they are isomorphic.
\hfill$\Box$

\

To classify the $q$-conjugacy classes in $\wtl{GL}_n(k((t)))$ is equivalent to compute the fibre of the following map $\pi:$
\[
q\rm{-conjugacy classes in }\wtl{GL}_n(k((t)))\stackrel{\pi}{\lra} q\rm{-conjugacy classes in }GL_n(k((t))).
\]
For simplicity assume that $k$ is algebraically closed. Our result is

\begin{thm}
Let ${\cal O}$ be a $q$-conjugacy class in $GL_n(k((t)))$ which corresponds to the $D_q$-module
$V=\bigoplus^l_{i=1}J_{d_i}\otimes V_{r_i,s_i}(\la_i)$, then

$\mathrm{(i)}$ $\pi^{-1}({\cal O})$ consists of a single class if $r_i\neq 0$ for some $i$,

$\mathrm{(ii)}$ $\pi^{-1}({\cal O})$ is a $k^\times$-torsor if all $r_i$'s are zero.
\end{thm}
\e{Proof.} Pick up $g\in{\cal O}.$ First note that the conjugate $(h,\om_h)(g,\om_g)\s_q (h,\om_h)^{-1}$ does not depend on the choice of $\om_h$, hence we may write it as $\wtl{\rm{Ad}}_qh(g,\om_g).$ If $\rm{Ad}_q h_1(g)= \rm{Ad}_q h_2(g),$ then $h:=h_2^{-1}h_1$ lies in the twisted centralizer of $g$, i.e. $\rm{Ad}_q h(g)=g.$ Equivalently $h$ corresponds to an automorphism of the $D_q$-module $V$.

If $r_i\neq 0,$ let $h$ corresponds to the automorphism of $V$ which acts by a non-zero scalar $\a$ on $J_{d_i}\otimes V_{r_i,s_i}(\la_i)$ and by identity on other summands, then $\wtl{\rm{Ad}}_q h(g,\om_g)=(g,\a^{-r_i d_i}\om_g).$ This proves (i) since $\a$ is arbitrary.

If all $r_i$'s are zero, then we may assume $g\in GL_n(k)$, $\om_g\in k^\times$. Then twisted centralizers of $g$ also lie in $GL_n(k),$ hence fix $(g,\om_g).$ This proves (ii).
\hfill $\Box$

\section{\large Conjugacy classes in symplectic and orthogonal groups}

In this section we will give our classification of symplectic
(resp. orthogonal) $D_q$-modules over $k((t))$, whose isomorphism
classes are in one-to-one correspondence with the $q$-conjugacy
classes in symplectic (orthogonal resp.) groups over $k((t)).$
Throughout this section we assume that $k$ is algebraically closed of characteristic zero, and $q\in k^\times$ is not
a root of unity.

\begin{Def}
A symmplectic $($resp. orthogonal$)$ $D_q$-module over $k((t))$ is a triple $(V,\Phi,\lg,\rg)$ such that $(V,\Phi)\in{\cal D}_q(k,t)$ and $\lg,\rg$
is a non-degenerate symplectic $($resp. symmetric$)$ form on $V$ such that $\lg\Phi v,\Phi w\rg=\s_q\lg v,w\rg, \forall v, w\in V.$
\end{Def}

We can define the category ${\cal D}_q(k,t)_S$ (resp. ${\cal D}_q(k,t)_O$) of symplectic (resp. orthogonal) $D_q$-modules, where morphisms
 between two objects $(V,\Phi,\lg,\rg)$ and $(V',\Phi',\lg,\rg ')$ is given by $k((t))$-linear maps $\rho: V\ra V'$ such that
$\Phi'\circ\rho=\rho\circ\Phi$ and $\lg \rho v,\rho w\rg '=\lg v,w\rg.$

Given an $2n$(resp. $n$)-dimensional vector space $V$ over $k$ and a non-degenerate symplectic (resp. symmetric) bilinear form $\lg,\rg$ on $V$, there is a unique extension of $\lg,\rg$ to $V_{k((t))}:=V\ot_k k((t))$ by linearity. We define the group $Sp_{2n}(k((t)))$ (resp. $O_n(k((t)))$) to be the $k((t))$-points of $Sp(V)$ (resp. $O(V)$). Then $\s_q$ acts on this group by $(\s_q g)(v\ot \a):=\a(\rm{id}_V\ot \s_q)(gv)$ for all $v\in V, \a\in k((t)).$

Let ${\cal D}_q^0(k,t)_S$ (resp. ${\cal D}_q^0(k,t)_O)$ be the full subcategory of ${\cal D}(q,t)_S$ (resp. ${\cal D}(q,t)_O$), which consists of $D_q$-modules with bilinear form arising in this way, i.e. arising from scalar extension of a bilinear form on a vector space over $k$.

The notions defined above are independent of the choice of $\lg,\rg$ on $V$, since we assume that $k$ is algebraically closed.

\begin{prop}\label{prop}
There is a bijection between the set of isomorphism classes of $2n$(resp. $n$)-dimensional $D_q$-modules in
${\cal D}^0_q(k,t)_{S}$ $($resp. ${\cal D}^0_q(q,t)_O)$ and the set of $q$-conjugacy classes in $Sp_{2n}(k((t)))$
$($resp. $O_n(k((t))))$.
\end{prop}
\e{Proof.} For $g\in Sp_{2n}(k((t)))$ we have constructed the $D_q$-module $(V_{k((t))},\Phi_g)$ in the proof of Proposition \ref{1}. But we have $\lg \Phi_g v,\Phi_g w\rg= \lg g\circ(\rm{id}_V\ot \s_q) v,
 g\circ(\rm{id}_V\ot\s_q) w\rg=\lg(\rm{id}_V\ot\s_q)v,(\rm{id}_V\ot\s_q)w\rg=\s_q\lg v,w\rg$ for all $v,w\in V_{k((t))}.$ This implies that $(V_{k((t))},\Phi_g,\lg,\rg)\in {\cal D}^0_q(k,t)_S.$ Conversely, given a $D_q$-module $(W,\Phi,\lg,\rg)$ in ${\cal D}^0_q(k,t)_S$, by assumption we can write $W=V_{k((t))}$ for some $k$-vector space $V$ and $\lg,\rg$ is obtained from scalar extension of a symplectic form on $V.$ Then $\lg \Phi v,\Phi w\rg=\s_q\lg v,w\rg=\lg v,w\rg$ for any $v,w\in V$, since $\lg v,w\rg\in k.$ Therefore we get an element in $Sp(V_{k((t))})=Sp_{2n}(k((t))).$ It is easy to check that these constructions give the required bijection. The proof for the orthogonal case is similar.
\hfill$\Box$

\

Let us denote $k^\times /q^\mb{Z}$ by $E_q(k),$ and the $2$-torsion points of $E_q(k)$ by $E_q[2]$. Then $E_q[2]\h (\mb{Z}/2\mb{Z})^2.$ In fact $E_q[2]=\{\pm 1 \mod q^\mb{Z}, \pm\sqrt{q}\mod q^\mb{Z}\}$. If $k=\mb{C}$ or $\overline{\mb{Q}}_p$ for some prime $p$, and $q\in k^\times $ with $|q|<1,$ then $E_q(k)$ can be identified with the $k$-points of Tate's elliptic curve $E_q=\mb{G}_m/q^\mb{Z}$ defined
over $\mb{Z}[[q]],$ see \cite{T}.

The main result of this section is the following

\begin{thm}\label{thm}
Any symplectic $($resp. orthogonal$)$ $D_q$-module has an orthogonal decomposition such that each summand is isomorphic to one of the following objects.

$\mathrm{(i)}$ $J_d\otimes V_{r,s}(\la)\oplus J_d\otimes V_{-r,s}(\la^{-1})$, where either $\bar{\la}\not\in E_q[2]$ or $r\neq 0$, $J_d\otimes V_{r,s}(\la)$ and
$J_d\otimes V_{r,s}(\la^{-1})$ are isotropic and dual to each other with respect to $\lg,\rg.$

$\mathrm{(ii)}$ $J_d\otimes V_{0,1}(\la)\oplus J_d\otimes V_{0,1}(\la)$ where $\bar{\la}\in E_q[2]$ and $d$ is odd $($resp. even$)$, and the two copies of
$J_d\otimes V_{0,1}(\la)$ are isotropic and dual to each other with respect to $\lg,\rg.$

$\mathrm{(iii)}$ $J_d\otimes V_{0,1}(\la)$ where $\bar{\la}\in E_q[2]$ and $d$ is even $($resp. odd$)$.
\end{thm}

\begin{Remark}
Proposition \ref{prop} and Theorem \ref{thm} imply that the natural map
\[
q\rm{-conjugacy classes in }Sp_{2n}(k((t)))\lra q\rm{-conjugacy classes in }GL_{2n}(k((t)))
\]
is injective. The same is true for the orthogonal group.
\end{Remark}

The rest of this section is devoted to the proof of Theorem \ref{thm}. We start from a lemma on the decomposition of $D_q$-modules.

\begin{lemma}\label{pair}
If $V=V_1\op V_2$ is a decomposition of $D_q$-module with $V_1$ indecomposable and isotropic, then there exists a $D_q$-submodule $V_1'$ of $V$ such
that $V_1\cap V_1'=0$ and $V_1'\h V_1^\vee:=\Hom_{k((t))}(V_1,k((t)))$ via the pairing $\lg,\rg.$ In particular $V_1'$ is also indecomposable and isotropic,
and $V_1\op V_1'$ is non-degenerate.
\end{lemma}
\e{Proof.} Let $v$ be a cyclic vector of $V_1$, then $v,\Phi
v,\ld,\Phi^{d-1}v$ is an $k((t))$-basis of $V_1$, where $d=\dim V_1.$
$V\h V^\vee$ via $\lg,\rg$,
 hence we can pick up $v'\in V$ such that $\lg v',v\rg=1,$ $\lg v',\Phi^i v\rg=0,$ $i=1,\ld,d-1,$ and $v'\perp V_2.$ Let $V_1'$ be the $D_q$-module
 generated by $v'.$

We first prove that $\dim V_1'=d.$ $V_1'\ss V_2^\perp $ implies that
$\dim V_1'\leq d.$ For the other direction we prove that $v',\Phi
v',\ld,\Phi^{d-1}v'$
 are linearly independent over $k((t))$. Suppose that
\[
a_0 v'+a_1\Phi v'+\cs+a_{d-1}\Phi^{d-1}v'=0,
\]
where $a_i\in k((t)).$ Paring with $\Phi^{d-1}v,\ld,v$ successively, we
see that $a_{d-1}=\cs=a_0=0.$ Hence $\dim V_1'\geq d.$

The map $V_1'\ra V_1^\vee, v\mt\l\lg v,\c\rg\r|_{V_1}$ is injective, hence an isomorphism. $V_1\cap V_1'=\emptyset$ since $V_1$ is
isotropic. The rest of the lemma are easy to verify.
\hfill $\Box$

\

Let $V$ decompose as in Corollary \ref{cor}. Namely, assume that
 \[
 V\h \bigoplus^l_{i=1}J_{d_i}\otimes V_{r_i,s_i}(\la_i).
 \]
 We shall show that most of the indecomposable summands are isotropic and orthogonal to each other.

Recall that we have a functor
\[
ext_n:\q {\cal D}_q(k,t)\longrightarrow{\cal D}_{q^{1/n}}(k,t^{1/n}).
\]
Let us consider $ext_s(J_d\otimes V_{r,s}(\la))$, where $J_d\otimes V_{r,s}(\la)$ is a summand of $V$. It is easy to check that
  \[
  ext_s(V_{r,s}(\la)) \h \bigoplus^s_{i=1}V_i,
    \]
where $V_i=k((t^{1/s}))v_i\in {\cal D}_{q^{1/s}}(k,t^{1/s})$ is one-dimensional such that $\Phi v_i=\la_i t^{r/s} v_i,$  and $\la_1,\ld,\la_s$ are distinct
  $s$-th roots of $\la$. Let us still write $J_d$ for $ext_s(J_d)$. Then
 \[
ext_s(J_d\otimes V_{r,s}(\la))=\bigoplus^s_{i=1} J_d\otimes v_i.
 \]

 \begin{lemma}\label{lemma}
 If we have another summand $J_{d'}\otimes V_{r's'}(\la')$ of $V$ with decomposition
 \[
ext_{s'}(J_{d'}\otimes V_{r',s'}(\la'))=\bigoplus^{s'}_{j=1} J_{d'}\otimes v'_j,
 \]
then $J_d\otimes v_i\perp J_{d'}\otimes v'_j$ unless $s=s', r=-r'$ and $\la_i\la'_j\in q^{\mb{Z}/s}.$
\end{lemma}

Before we prove this lemma let us observe an elementary fact:

\begin{lemma}\label{fact}
If $A\in M_{d\times d'}(k((t)))$, $\la\in k,$ $\nu\in\mb{Z}$, $\la t^\nu J_d A J_{d'}^t=\s_q A,$ then $A=0$ unless $\nu=0$ and $\la=q^\xi$ for some
$\xi\in\mb{Z},$ in which case $A\in t^\xi M_{d\times d'}(k).$ Here $J_d$ $($abuse of notation$)$ denotes the $d\times d$ Jordan block.
\end{lemma}
\e{Proof.}
Write $A=\sum^\i_{i=i_0}A_i t^i$, $A_i\in M_{d\times d'}(k),$ $A_{i_0}\neq0,$ then $J_d A_{i-\nu}J_{d'}^t=q^i A_i$, $i\in\mb{Z}.$ Suppose $\nu\neq0$,
say $\nu>0,$ then one finds that $A_{i_0}=0,$ a contradiction. Now suppose $\nu=0.$ Then $\la J_d A_i J_{d'}^t=q^iA_i.$

In general for $U\in M_{d\times d}(k), V\in M_{d'\times d'}(k)$, the
linear transform $\phi_{UV}: M_{m\times m'}(k)\ra M_{d\times
d'}(k),$ $A\mt UAV$ has eigenvalues contained in $Sp(U)Sp(V):
=\{\la\mu| \la\in Sp(U), \mu\in Sp(V)\}.$ To prove this we
may restrict to the case that $U, V$ lie in the open dense (Zariski topology) subset of semisimple
elements of $M_{d\times d}(k)$. Note that if $uU^t =\la u,$ $vV=\mu
v$, then $\phi_{UV}(u^t v)=\la\mu u^tv.$ Apply this fact to $J_d$
and $J_{d'}^t$ one sees that the claim follows. \hfill $\Box$

\

{\flushleft{\e{Proof.}}} (of Lemma \ref{lemma})

Let $\{u_1,\ld,u_d\}$ be a basis of $J_d$ such that $\Phi|_{J_d}$ is
given by the matrix $J_d.$ Similarly we choose the basis for
$J_{d'}.$ Let $A$ be the $d\times d'$ matrix $(\lg u_k\ot v_i,
u'_l\ot v_j'\rg)_{k,l}\in M_{d\times d'}(k((t^{1/[s,s']})))$.
Applying $\phi$ gives
\[
\la_i\la_j't^{r/s+r'/s'}J_dAJ_{d'}^t=\s_q A.
\]
If $A\neq 0$, then previous lemma implies that $s=s', r=-r'$, and $\la_i\la_j'=q^{\xi/s}$ for some $\xi\in\mb{Z},$ and $A\in t^{\xi/s}M_{d\times d'}(k)$.
\hfill $\Box$

\

{\flushleft{\e{Proof.}}} (of Theorem \ref{thm})

Lemma \ref{lemma} implies that $J_d\otimes V_{r,s}(\la)$ is isotropic unless $r=0, s=1$ and $\bar{\la}\in E_q[2],$ and
$J_d\otimes V_{r,s}(\la)\perp J_{d'}\otimes V_{r',s'}(\la')$ unless $s=s',$ $r=-r'$ and $\la\la'\in q^\mb{Z}.$ A moment's
thought gives that there is a decomposition of the form
\[
V=\l(\bigoplus_{r\neq 0\rm{ or }\bar{\la}\not\in E_q[2],d,n}J_d\otimes V_{r,s}(\la)\r)\bigoplus\bigoplus_{\bar{\la}\in E_q[2]}
\l(\bigoplus_dJ_d\otimes V_{0,1}(\la)\right).
\]
such that the outmost summands are mutually orthogonal. Apply Lemma \ref{pair} and Lemma \ref{lemma} to the first summand repeatedly
we obtain case (i) of Theorem \ref{thm}.

Next we treat the remaining parts. For each $\la\in E_q[2]$, by previous constructions and lemmas we know that there exists a basis $\{v_\a\}$
of \[
\bigoplus_dJ_d\otimes V_{0,1}(\la)
\]
such that

 (a) $\Phi$ acts on these basis by a constant matrix, whose semisimple part is the scalar $\la$,

(b) the pairing $\lg,\rg$ on these basis takes values in $t^\xi k$.
\\
Then $\{t^{-\xi/2}v_\a\}$ form a basis of
$ext_2\l(\bigoplus\limits_d J_d\otimes V_{0,1}(\la)\r)$, and the pairing on the basis takes values in $k$. Let $V$ be the $k$-span of $\{t^{-\xi/2}v_\a\}$,
then $\Phi$ acts on $V$ as a unipotent map, note that $\la^2=q^\xi.$

Theory of nilpotent orbits \cite{CM} implies that
$V=\bigoplus_i V_i$ an orthogonal decomposition such that each $V_i$
is one of the following, according to $\lg,\rg$ symplectic (resp.
symmetric):

(1) an indecomposable $\Phi$-module of even (resp. odd) dimension,

(2) a direct sum of two copies of isotropic indecomposable $\Phi$-module of odd (resp. even) dimension, which are dual to each other.

It is clear that
\[
\bigoplus_d J_d\otimes V_{0,1}(\la)=\bigoplus_i(t^{\xi/2}V_i)\otimes_k k((t))
\]
gives the decomposition where the summands occur as case (ii) or (iii) of Theorem \ref{thm}.
\hfill $\Box$

\section{\large Principal bundles over elliptic curves}

We shall discuss the connection with holomorphic principal bundles over elliptic curves. In \cite{vdP-R} it is observed that for $q\in\mb{C}^\times, |q|<1,$ there is a functor from the category of split $D_q$-modules (split means $V\h gr V$) over $\mb{C}\{t\}$, the field of convergent Laurent series, to the category of vector bundles on $E_q(\mb{C})\h \mb{C}^\times /q^\mb{Z}$, which is bijective on isomorphism classes of objects and respects tensor product. Moreover it is also proved that split $D_q$-modules over $\mb{C}\{t\}$ have the same classification as $D_q$-modules over $\mb{C}((t)).$

In summary, there is a formal correspondence between isomorphism classes of $D_q$-modules over $\mb{C}((t))$ and vector bundles over $E_q.$ Now let us explain the relations with Atiyah's theorem \cite{A} and the work of Friedman and Morgan \cite{FM}.

The $D_q$-module $J_d$ corresponds to the unique indecomposable vector bundle $I_d$ of rank $n$ over $E_q$ all of whose Jordan-Holder constituents are isomorphic to ${\cal O}_{E_q}.$

If we denote by $M_n(r,\la)$ the $D_q$-module $\mb{C}((t))_q[T]/(T^n+(-1)^{n-r}\la^{-1}t^r)$, then $M_n(r,\la)$ corresponds to the unique indecomposable bundle $L_n(r,\bar{\la})$ over $E_q$ of rank $n$ such that $\det L_n(r,\bar{\la})={\cal O}_{E_q}((r-1)p_{\bar{1}}+p_{\bar{\la}})$, where $p_{\bar{1}}$ and $p_{\bar{\la}}$ are divisors on $E_q$ corresponding to $\bar{1}$ and $\bar{\la}$ in $ \mb{C}^\times/q^\mb{Z}$ respectively.

Atiyah's theorem states that $I_d\ot L_n(r,\bar{\la})$ exhaust all indecomposable vector bundles over $E_q$ up to isomorphism, and any vector bundle over $E_q$ splits into a direct sum of indecomposables ones, with the quadruples $(d,n,r,\bar{\la})$ uniquely determined.

Thus the notions of degree, rank, stability and Harder-Narasimhan filtration on both categories coincide. For instance, any semi-stable vector bundle of slope $r/n$ on $E_q$ is of the form $\bigoplus_i I_{d_i}\ot L_n(r,\bar{\la}_i).$ Compare this with Section 2.

 Inspired by the notion of regular vector bundles defined in \cite{FM}, let us make the following
 \begin{Def}
 For a semi-stable $V\in {\cal D}_q(\mb{C},t)$, ${\cal D}_q(\mb{C},t)_S$ or ${\cal D}_q(\mb{C},t)_O$, define
 \[
 \wtl{gr}V:=\bigoplus_i V_i/V_{i-1}
 \]
 where $V_i\subseteq V$ form a filtration of $D_q$-submodules such that $V_i/V_{i-1}$ are stable with $\mu(V_i/V_{i-1})=\mu(V).$ $V_1$ and $V_2$ are called S-equivalent if $\wtl{gr}V_1\h \wtl{gr}V_2.$ $V$ is said to be regular if $\dim \mathrm{Aut}V\leq \dim\mathrm{Aut}V'$ for any $V'$ which is $S$-equivalent to $V$.
 \end{Def}

Of course we take the automorphism group in the corresponding categories. Lemma \ref{aut} implies that $\rm{End}(V_{r,s}(\la))\h\mb{C},$ $\rm{End}(J_d\ot V_{r,s}(\la))\h \mb{C}[t]/(t^d)$ in ${\cal D}_q(\mb{C},t)$.
Let us give some results on regular semi-stable $D_q$-modules, which are compatible with the description of regular bundles in \cite{FM}. It is easy to figure out the proof, which we shall not give here.

\begin{thm}
A semi-stable $V\in {\cal D}_q(\mb{C},t)$  is regular if and only if it is isomorphic to
\[
\bigoplus_i J_{d_i}\ot V_{r,s}(\la_i)
\]
where $\bar{\la}_i$'s  are all distinct. In this case $\dim \mathrm{Aut} V=\sum_i d_i.$
\end{thm}

\begin{thm}
A semi-stable $V\in {\cal D}_q(\mb{C},t)_{S}$ is regular if and only if the image of $V$ in ${\cal D}_q(\mb{C},t)$ is regular, i.e. isomorphic to
\[
\bigoplus_{i, \bar{\la}_i\not\in E_q[2]}J_{d_i}\ot V_{0,1}(\la_i)\oplus J_{d_i}\ot V_{0,1}(\la_i^{-1})\bigoplus \bigoplus_{\bar{\la}\in E_q[2]}J_{2d_{\bar{\la}}}\ot V_{0,1}(\la),
\]
where $(\bar{\la}_i,\bar{\la}_i^{-1})$'s are all distinct. In this case
\[
\dim\mathrm{Aut}V=\sum_i d_i+\sum_{\bar{\la}\in E_q[2]}d_{\bar{\la}}=\f{1}{2}\dim V.
\]
\end{thm}

We refer the reader to \cite{FM} Proposition 7.7 for the orthogonal case, which is more complicated.

\

Finally we would like to mention that, if we consider a general complex connected semisimple algebraic group $G$, then along this line the theorem due to of V. Baranovsky and V. Ginzburg \cite{BG}, states that there is a bijection between the set of integral $q$-conjugacy classes (i.e. containing an element of $G(\mb{C}[[t]])$) in $G(\mb{C}((t)))$
and the set of isomorphism classes of semi-stable holomorphic principal $G$-bundles on $E_q.$

Address:  Dept. of Math., Hong Kong University of Science and Technology.

Email:  ldxab@ust.hk

\end{document}